\documentclass{amsproc}

\newtheorem{Theorem}{Theorem}
\newtheorem{Corollary}[Theorem]{Corollary}
\newtheorem{theorem}{Theorem}[section]
\newtheorem{lemma}[theorem]{Lemma}

\theoremstyle{definition}

\newtheorem{assumption}[theorem]{Assumption}

\theoremstyle{remark}

\numberwithin{equation}{section}

\newcommand{\R}{{\mathbb R}}
\newcommand{\Z}{{\mathbb Z}}
\newcommand{\FF}{{\mathcal F}}

\newcommand{\Int}{{\rm Int}}
\begin{document}

\title[Prime end rotation numbers 
of invariant separating continua]{Prime end rotation numbers 
of invariant separating continua of annular homeomorphisms}

\author{Shigenori Matsumoto}
\address{Department of Mathematics, College of
Science and Technology, Nihon University, 1-8-14 Kanda, Surugadai,
Chiyoda-ku, Tokyo, 101-8308 Japan
}
\email{matsumo@math.cst.nihon-u.ac.jp
}
\thanks{The author is partially supported by Grant-in-Aid for
Scientific Research (C) No.\ 20540096.}
\subjclass{Primary 37E30,
secondary 37E45.}

\keywords{continuum, rotation set, prime end rotation number,
Brouwer line, foliations}

\date{\today }

\begin{abstract}
Let $f$ be a homeomorphism of the closed annulus $A$ isotopic to the identity,
and let $X\subset {\rm Int}A$ be an $f$-invariant
 continuum which separates $A$ into
two domains, the upper domain $U_+$ and the lower domain $U_-$.
Fixing a lift of $f$ to the universal cover of $A$,
one defines the rotation set $\tilde \rho(X)$ 
of $X$ by means of the invariant probabilities
on $X$, as well as the prime end rotation number
$\check\rho_\pm$ of $U_\pm$.
The purpose of this paper is to show that $\check\rho_\pm$
belongs to $\tilde\rho(X)$ for any separating invariant
continuum $X$.

\end{abstract}

\maketitle

\section{Introduction}
Let $f$ be a homeomorphism of the closed annulus $A=S^1\times[-1,1]$,
isotopic to the identity, i.\ e.\ $f$ preserves the orientation and 
each of the boundary
 components
$\partial_\pm A=S^1\times\{\pm 1\}$. Suppose there is an $f$-invariant partition
of $A$; $A=U_-\cup X \cup U_+$, where $U_{\pm}$ is a connected
open set containing the boundary component $\partial_\pm A$
and $X$ is a connected compact set. Let 
$$
\pi: \tilde A=\R\times [-1,1]\to S^1\times[-1,1]
$$
be the universal covering map and $T:\tilde A\to \tilde A$
a generator of the covering transformation group;
$T(\xi,\eta)=(\xi+1,\eta)$. Denote by $p:\tilde A\to \R$ the projection onto
the first factor.

Fix once and for all  a lift $\tilde f:\tilde A\to\tilde A$ of $f$.
Then the function 
$p\circ\tilde f-p$ is $T$-invariant and can be looked upon as a function
on the annulus $A$.
Define the {\em rotation set} $\tilde \rho(X)$ as
the set of values $\mu(p\circ\tilde f-p)$, where $\mu$ ranges over
the $f$-invariant probability measures supported on $X$.
The rotation set is a compact interval (maybe one point) in $\R$, which
depends upon the choice of the lift $\tilde f$ of $f$.

The first example of an invariant continuum $X$ such that 
the frontiers of $U_\pm$ satisfy
${\rm Fr}(U_+)={\rm Fr}(U_-)=X$ 
and that the rotation set $\tilde\rho(X)$ is not a
singleton is constructed by G. D. Birkhoff in his 1932 year paper
\cite{B}, and is refered to as a {\em Birkhoff attractor}.
It turns out that the 
Birkhoff attractor is an indecomposable continuum (\cite{C,L2}).
Furthermore it is shown by P. Le Calvez (\cite{L1})
that for any rational number between the two prime end
rotation numbers is realized by a correspoding periodic 
point of $\tilde f$.

Let $\hat U_\pm=U_\pm\cup\partial_\infty U_\pm$ 
be the prime end compactification of $U_\pm$,
where $\partial_\infty U_\pm$ is the space of the prime ends (\cite{E,M,MN}). 
The space $\partial_\infty U_\pm$ is homeomorphic to the circle and
$\hat U_\pm$ to the closed annulus. As is well known,
the homeomorphism $f$ restricted to  $U_\pm$ extends to a homeomorphism
$\hat f_\pm:\hat U_\pm\to\hat U_\pm$.
Denoting $I_+=[0,1]$ and $I_-=[-1,0]$, define a homeomorphism
$$
\Psi_{\pm}:\hat U_\pm\to S^1\times I_\pm$$
such that
$\Psi_\pm(\partial_\infty U_\pm)=S^1\times 0$.
By some abuse of notations denote by
$\pi:\check U_\pm\to \hat U_\pm$
the universal covering map. Thus $\pi^{-1}(U_\pm)$ is considered
to be a subspace of both
$\tilde A$ and $\check U_\pm$.   
Let $\check f_\pm:\check U_\pm\to \check U_\pm$ be the lift of 
$\hat f_\pm$
such that $\check f_\pm=\tilde f$ on $\pi^{-1}(U_\pm)$. 
The rotation number of the restriction of $\check f_\pm$ to
$\pi^{-1}(\partial_\infty U_\pm)$, denoted by $\check\rho_\pm$,
is called the {\em prime end rotation number} of $U_\pm$.

The purpose of this paper is to show the following.

\begin{Theorem} \label{tx}
The prime end rotation number $\check\rho_\pm$ belongs to $\tilde\rho(X)$.
\end{Theorem}

This result is already known for $X={\rm Fr}(U_-)={\rm Fr}(U_+)$
(\cite{BG}), and for any $X$ if the homeomorphism $f$ is area
preserving (Lemma 5.4, \cite{FL}).

It is shown in Theorem 2.2 of \cite{F} 
that any rational number in $\tilde\rho(X)$
is realized by a periodic point if $X$ consists of nonwandering points.
Notice that then $X$, consisting of chain recurrent points, is chain
transitive since it is connected, and thus satisfies the condition
of Theorem 2.2.
As a corollary we have

\begin{Corollary} \label{cx}
If $X$ consists of nonwandering points and if $p/q$ 
lies in the closed interval bounded by $\check\rho_-$ and
 $\check\rho_+$, then there is a point $x\in\pi^{-1}(X)$ such that
 $\tilde f^q(x)=T^p(x)$.
\end{Corollary}

In what follows we also use the following notation. Let
$$\check\Psi_\pm:\check U_\pm\to \R\times I_\pm$$ 
be a lift of $\Psi_\pm$,
and define $\check p_\pm:\check U_\pm\to\R$ by 
$\check p_\pm=p\circ\check\Psi_\pm$. 
The projection $\check p_\pm$ is within a bounded
error of $p$ on $\pi^{-1}(C)$ for a compact domain $C$ of $U_\pm$.
But they may be quite different on the whole $\pi^{-1}(U_\pm)$.

\section{Proof}

First of all let us state a deep and quite useful 
theorem of P. Le Calvez (\cite{L3})
which plays a key role in the proof.
A fixed point free and orientation preserving homeomorphism $F$ of the
plane $\R^2$ is called a {\em Brouwer homeomorphism.} A proper oriented
simple curve $\gamma:\R\to\R^2$ is called a {\em Brouwer line} for
$F$ if $F(\gamma )\subset R(\gamma)$ and 
$F^{-1}(\gamma)\subset L(\gamma)$,
where $R(\gamma)$ (resp.\ $L(\gamma)$) is the right (resp.\ left) side complementary
domain
of $\gamma$, which is decided by the orientation of $\gamma$.

\begin{theorem} \label{t3}
Let $F$ be a Brouwer homeomorphism commuting with the elements of a
 group
$\Gamma$ which acts on $\R^2$ freely and properly discontinuously.
Then there is a $\Gamma$-invariant oriented topological foliation of $\R^2$
whose leaves are Brouwer lines of $F$.
\end{theorem}

The proof of Theorem \ref{tx} is by absurdity. 
Assume in way of contradiction that
$\check\rho_-<p/q<\inf\tilde\rho(X)$.
Considerng $\tilde f^qT^{-p}$ instead of $\tilde f$, it suffices to
deduce a contradiction under the following assumption.

\begin{assumption} \label{ax}
$\check\rho_-<0<\inf\tilde\rho(X)$.
\end{assumption}

Since $\inf\tilde\rho(X)>0$, the map 
$\tilde f$ does not admit a fixed point in $\pi^{-1}(X)$.
The overall strategy of the proof is to modify the homeomorphism
$f$ away from $X$ to a new one $g$
without creating fixed points in $A$ such  that
the restrictions of $\tilde g$ to the 
lifts of the both boundary circles $\pi^{-1}(\partial_\pm A)$
are nontrivial rigid translations by the same translation number.
Then by glueing the two boundary circles we obtain a torus
$T^2$ and a homeomorphism on $T^2$. Now we can apply Theorem
\ref{t3} to the lift of the homeomorphism to the universal covering
space. This yields a topological foliation on $T^2$, which 
has long been well understood. The proof will be done
by analyzing the foliation.
We first prepare a lemma which is necessary for the desired
modification. We do not presume Assumption \ref{ax} in the
following.

\begin{lemma} \label{l1}
Assume $\tilde f$ does not admit a fixed point in 
$\pi^{-1}(X)$.
Then the prime end rotation number $\check\rho_\pm$ is nonzero.
\end{lemma}

{\sc Proof:} 
Consider the mapping $\tilde f-{\rm Id}$ defined on $\tilde A$.
Since it is $T$-invariant, it yields a mapping from $A$, still denoted
by the same letter. Then since there is no fixed point of $\tilde f$
in $X$, we have $(\tilde f-{\rm Id})(X)\subset\R^2\setminus\{0\}$.
Therefore there is an annular open neighbourhood $V$ of $X$ for which we
get a mapping
$$
\tilde f-{\rm Id}:V\rightarrow \R^2\setminus\{0\}.$$
Clearly for any positively oriented essential simple closed curve $\gamma$ in
$V$, the degree of the map
$$
\tilde f-{\rm Id}:\gamma\rightarrow \R^2\setminus\{0\}$$
must be the same.
If the curve $\gamma$ is contained in $U_{\pm}$, then the
degree can be studied by considering the map $\check f_\pm$
defined on the lift $\check U_\pm$ of the prime end
compactification $\hat U_\pm$. If the prime end
rotation number $\check\rho_\pm$ is nonzero, the
degree is clearly 0. Notice that our definition of the degree
differs from the usual definition of the index.

To analyze the case $\check\rho_\pm=0$,
we need the following form of
the Cartwright-Littlewood theorem \cite{CL}.

\begin{theorem} \label{t4}
If $\check\rho_+=0$ and if ${\rm Fix}(\tilde f)\cap \pi^{-1}(X)=\emptyset$,
then the map $\hat f_+$ on $\partial_\infty U_+$ is Morse Smale
and the attractors (resp.\ repellors) of $\hat f_+\vert_{\partial_\infty U_+}$
are attractors (resp.\ repellors) of the whole map $\hat f_+$.
\end{theorem}

This is slightly stronger than the usual version in which it is
assumed that Fix$(f)\cap X=\emptyset$. However the proof works as well
under the assumption of Theorem \ref{t4}. See e.\ g.\ Sect.\ 3 
of \cite{MN}.

Let us complete the proof of Lemma \ref{l1}. Theorem \ref{t4}
enables us to compute
the degree of the curve $\delta$ in $U_\pm$ when $\check\rho_\pm=0$. 
The degree 
is $n$ if $\delta\subset U_-$ and $-n$ if $\delta\subset U_+$,
where $n$ is the number of the attractors.
Since the degree must be the same in $U_-$ and $U_+$,
the conclusion follows.
\qed 
\smallskip

Now we have $\check\rho_-<0$ and $\check\rho_+\neq0$ by
Assumption \ref{ax} and Lemma \ref{l1}. 
Let us start the modification of $f$. 

\begin{lemma} \label{l2}
Under Assumption \ref{ax},
there exists a homeomorphism $g$ of $A$ such that
\\
(1) $g=f$ in some neighbourhood of $X$, 
\\
(2) $\tilde g$ does not admit a fixed point in $\tilde A$,
where $\tilde g$ is the lift of $g$ such that
$\tilde g=\tilde f$ on $\pi^{-1}(X)$, 
\\
(3) $\tilde g$ is a negative rigid translation by the
same translation number on $\pi^{-1}(\partial_\pm A)$, and
\\
(4) $\check p_-\circ \check g_--\check p_-\leq -c$ on $\hat U_-$
for some positive number $c$.
\end{lemma}

{\sc Proof:}
The modification in $U_-$ will be done in the following
way.
We identify $\hat U_-$ with $S^1\times[-1,0]$ by the homeomorphism
$\Psi_-$ and the universal covering space $\check U_-$ with
$\R\times[-1,0]$. Thus $\check p_-$ is just the projection onto
the first factor;
$\check p_-(\xi,\eta)=\xi$.
Since $\check\rho_-<0$, the lift 
$$
\check f_-:\R\times [-1,0]\to\R\times[-1,0]
$$ of $\hat f_-$ satisfies that $\check p_-\circ\check f_-(\xi,0)<\xi-2c$ for
some $c>0$. Therefore
changing the coordinates of $[-1,0]$ if necessary, one may assume
that $\check p_-\circ\check f_-(\xi,\eta)\leq \xi-c$ if 
$(\xi,\eta)\in\R\times[-1/2,0]$.
Define a homeomorphism $h$ of $S^1\times[-1,0]$ by
$$
h(\xi,\eta)=(\xi+\varphi(\eta){\rm mod}1,\eta),$$
where
$\varphi:[-1,0]\to(-\infty, 0]$ is a continuous function such that
$\varphi([-1/2,0])=0$ and 
$$
\varphi(\eta)\leq-\sup\{(\check p_-\circ \check f_--\check p_-)(\xi,\eta)\mid
\xi\in S^1\}-c.
$$
Define $g=f\circ h$. Then its lift $\check g_-$ satisfies
$$
\check p_-\circ\check g_--\check p_-\leq -c$$
on $\check U_-=\R\times[-1,0]$.
 Clearly condition (3) for $\pi^{-1}(\partial_-A)$ can be established by
a further obvious modification.

Now to modify $f$ in $U_+$, we do the same thing as in $U_-$.
If the prime end rotation number $\check\rho_+$ is negative,
then with an auxiliary modification we are done. If it is positive
 insert a time one map
of the Reeb flow.
\qed
\smallskip

Consider the torus $T^2$ which is obtained from $A$ 
by glueing the two boundary 
curves
$\partial_-A$ and $\partial_+A$. Then the condition (3) above shows that
$g$ induces a homeomorphism of $T^2$, again denoted by $g$. The
universal cover of $T^2$ is $\R^2$ and $\tilde A=\R\times[-1,1]$
is a subset of $\R^2$. The lift $\tilde g:\tilde A\to\tilde A$ can
be extended uniquely to a lift $\tilde g:\R^2\to\R^2$ of $g:T^2\to T^2$.
The covering transformation group $\Gamma$ is isomorphic to $\Z^2$,
generated by the horizontal translation $T$ and the vertical translation
by 2, denoted by $S$. Since $\tilde g$ is a Brouwer homeomorphism
which commutes with $\Gamma$, there is a
$\Gamma$-invariant oriented
foliation  on $\R^2$ whose leaves are Brouwer lines for $\tilde g$.
This yields an oriented foliation $\FF$ on the torus $T^2$.
The proof is divided into several cases according to the topological
type of the foliation $\FF$. We are going to deduce a contradiction
in each case. But before going into detail we need
another lemma.

\begin{lemma} \label{lz}
For any $C>0$ there is $n>0$ such that $p\circ\tilde g^n-p\geq C$ on $X$.
\end{lemma}

{\sc Proof:}
If not, there would be a point $x_n\in X$ for any $n>0$
such that
$$(p\circ \tilde g^{n}-p)(x_n)
=\sum_{j=0}^{n-1}(p\circ\tilde g-p)(g^j(x_n))< C$$
for some $C>0$,
 and the averages of Dirac masses
$$\mu_n=\frac{1}{n}\sum_{j=0}^{n-1}g^j_*\delta_{x_n}
$$
would satisfy $\mu_n(p\circ\tilde g-p)<C/n$.
Therefore an accumulation point $\mu$ of $\mu_n$
would have the property that
$\mu(p\circ\tilde g-p)\leq 0$, contradicting the
assumption
$\inf\tilde\rho(X)>0$.
\qed

\bigskip

{\sc Case 1}. {\em The foliation $\FF$ does not admit a compact leaf}.
Then $\FF$ is conjugate either to a linear foliation or to a Denjoy
foliation, both of irrational slope. The lift $\tilde \FF$ of $\FF$ 
to the open annulus
$\R^2/\langle T\rangle$ is conjugate to a foliation by vertical lines.
The space of leaves of $\tilde\FF$ is homeomorphic to $S^1$
and there is a projection 
from $\R^2/\langle T\rangle$ to $S^1$
along the leaves of the foliation. This lifts to a projection
$q:\R^2\to \R$. 

Now $q$ restricted to $\tilde A$
is within a bounded error of
the first factor projection
$p:\tilde A\to \R$ that we have used for the definition of 
the rotation set $\tilde\rho(X)$. In fact both $p$ and $q$ are
lifts of degree one maps from $\R^2/\langle T\rangle$ to $S^1$ and their
difference is bounded on the preimage $\tilde A=\pi^{-1}(A)$ of a compact
subset $A$. 
Thus Lemma \ref{lz} shows that
$q\circ\tilde g^n(x)\to\infty$ ($n\to\infty$) for $x\in\pi^{-1}(X)$.
That is, the foliation $\tilde\FF$ is oriented upward.
But this shows that $q\circ\tilde g(x)>q(x)$
even for a point $x\in\pi^{-1}(\partial_-A)$.
On the other hand  by condition (3) of Lemma \ref{l2},
$\tilde g$ is a negative translation on $\pi^{-1}(\partial_-A)$. 
A contradiction.

\smallskip

{\sc Case 2.1}. {\em The foliation $\FF$ admits a compact leaf $L$
of nonzero slope and does not admit a Reeb component}. 
In this case the lifted foliation $\tilde\FF$ is also
conjugate to the vertical foliation and the argument of Case 1 applies. 

\smallskip

{\sc Case 2.2}. {\em The foliation $\FF$ admits a Reeb component $R$
of nonzero slope}.
The Brouwer property of leaves implies that
$g(R)\subset{\Int}(R)$ or $g^{-1}(R)\subset{\Int (R)}$. 
That is, a point of the boundary of $R$ is
wandering under $g$. Therefore 
$\partial_-A$, 
consisting of nonwandering points of $g$ according to (3) of Lemma \ref{l2}, 
cannot intersect the boundary of $R$, which is however impossible since the
slope of $R$ is nonzero.
\smallskip

{\sc Case 2.3}. {\em The foliation $\FF$ admits a compact leaf of slope 0}.
Hereafter we only consider the dynamics and the foliation on the open
annulus $\R^2/\langle T\rangle$. Recall that $A$ is a subset of
$\R^2/\langle T\rangle$, and the homeomorphism $g$ on $A$ is extended
to the whole $\R^2/\langle T\rangle$, again denoted by $g$,
 in such a way that $g$ commutes
with the vertical translation $S$, while the foliation is denoted
by  $\tilde\FF$ as before. 

Now the foliation $\tilde\FF$ yields a partition $\mathcal P$ of the
open annulus $\R^2/\langle T\rangle$
into compact leaves, interiors of  Reeb components and  foliated $I$-bundles.
The set $\mathcal P$ is totally ordered by the height.
The minimal element which intersects $X$ cannot be a compact leaf
by the Brouwer line property. Let $R$ be the closure of the minimal
element. Thus $R$ is either a Reeb component or a foliated $I$-bundle
such that $\Int(R)\cap X\neq\emptyset$ and $\partial_-R\cap
X=\emptyset$,
where $\partial_-R$ is the lower boundary curve of $R$.

Assume for a while that $\partial_-R$ is oriented from the right
to the left. Thus the homeomorphism $g$ carries $\partial_-R$ into the
upper complement of $\partial_-R$.

\smallskip
{\sc Case 2.3.1} {\em $R$ is a Reeb component}. 
First notice that $g(R)\subset\Int R$ and that
the interior leaves of $R$ are oriented upwards
by the assumption $\inf\tilde\rho(X)>0$ and the fact that 
$g(X\cap R)\subset X\cap R$. 
Choose a simple arc 
$$\alpha:[0,1]\to \pi^{-1}(R)$$ such that
$\alpha(0)\in\pi^{-1}(\partial_-R)$,
$\alpha(1)=\tilde g(\alpha(0))$, and
$\alpha((0,1))\subset\Int(\pi^{-1}(R))\setminus\tilde g(\pi^{-1}(R))$.
Since $g^{-1}(\pi(\alpha))$ is below $\Int R$,  $\tilde g^{-1}(\alpha)$,
and hence $\alpha$, is contained in $\pi^{-1}(U_-)$.

Concatenating nonnegative iterates of $\alpha$, we obtain a simple path
$\gamma: [0,\infty)\to \pi^{-1}(R\cap U_-)$ such that $\tilde g\circ
\gamma(t)=\gamma(t+1)$
for any $t\geq 0$.  Let $q:\pi^{-1}(\Int(R))\to\R$ be the lift
of the projection along the leaves. 
Since $\gamma([1,\infty))$
is contained in the lift of a compact subset 
$\tilde g(R)\subset\Int(R)$ and the leaves in $\Int(R)$ is oriented
upward, we have
$q\circ\gamma(t)\to\infty$ as $t\to\infty$.
We also have $p\circ\gamma(t)\to\infty$ because $q$ is within bounded
error of $p$ on $\gamma([1,\infty))$.

On the other hand by condition (4) of Lemma \ref{l2}, we have
$\check p\circ\gamma(t)\to-\infty$ as $t\to\infty$.
In particular the curve $\gamma$ is proper both in $\tilde A$ and
in $\check U_-$ pointing toward the opposite direction. 
By joining the point $\gamma(0)$ to an appropriate point in 
$\pi^{-1}(\partial_-A)$, we obtain a simple curve $\delta$ in $\pi^{-1}(U_-)$
starting at a point on $\pi^{-1}(\partial_-A)$
which  extends $\gamma$.

Notice that there is a point of $\pi^{-1}(X)$ on the left of
a proper oriented curve $\delta$ in $\tilde A$,
because the map $p$ is bounded from below on $\delta$
and a high iterate of $T^{-1}$ carries a point in $\pi^{-1}(X)$
beyond that bound. (There might be a point of
$\pi^{-1}(X)$
on the right of $\delta$ however.)

Let $x$ be a point in
$\pi^{-1}(\partial_-A)$ left to the initial point of $\delta$.
Then there is a simple path
$\beta:[0,\infty)\to\pi^{-1}(U_-)$ such that $\beta(0)=x$,
$\lim_{t\to\infty}\beta(t)\in\pi^{-1}(X)$, and  $\beta$ is disjoint
from $\delta$. 
The path $\beta$, extendable in $\pi^{-1}(A)$ is also
extendable in $\check U_-$, the lift of the prime end
compactification. (See e.\ g.\ Lemma 2.5 of \cite{MN}.)
This implies that $\beta$ defines a simple path in
$\check U_-$ joining $x$ to a prime end in $\pi^{-1}(\partial_\infty
U_-)$ 
without intersecting
$\delta$, 
which is impossible since $\pi^{-1}(\partial_\infty U_-)$ 
is contained in the right side
of the proper path $\delta$ in $\check U_-$
since $\check p_-\delta (t)\to-\infty$,
while $x$ is on the left side.
A contradiction.

\smallskip
{\sc Case 2.3.2} {\em $R$ is a foliated $I$-bundle}.
Thus the upper boundary curve $\partial_+R$ of $R$ is also
oriented from the right to the left, and its image by $g$ lies
on the upper complement of $R$. The interior leaves of $R$ are
oriented upward.

Recall that the boundary component $\partial_-A$ consisting
of nonwandering points cannot intersect
a compact leaf. Moreover $\partial_-A$ lies in a Reeb component
or a foliated $I$-bundle  whose interior leaves are oriented
downward since $p\tilde g^n(x)\to-\infty$ as $n\to \infty$
for $x\in\pi^{-1}(\partial_-A)$.
Let $C$ be the annulus in $\R^2/\langle T \rangle$ bounded by
$\partial_-A$ and $\partial_+R$, the upper boundary curve of $R$.
Notice that $\Int(C)$ contains $\partial_-R$.

\smallskip
{\sc Case 2.3.2.1} {\em The intersection $X\cap C$ has a component which
separates $\partial_-A$ from $\partial_+A$}.
One can derive a contradiction by the same argument as in Case 2.3.1, since
the like defined path $\gamma$ cannot evade $R$.

\smallskip
{\sc Case 2.3.2.2} {\em There is a simple path in $U_-$ joining
a point in $\partial_-A$ with a point in $\partial_+R$.}
Notice first of all that $g^{-1}(C)\subset C$.
Let
$
{\mathcal Y}$ 
be the family of the connected components of $\pi^{-1}(X\cap C)$. Then any 
element $Y\in\mathcal Y$ is compact, and intersects $\pi^{-1}(\partial_+R)$
since otherwise $Y$ would be a connected component of $\pi^{-1}(X)$ itself.

Choose a simple curve $\gamma:[0,1]\to \pi^{-1}(C)$ such that
\\
(1) $\gamma(0)\in \pi^{-1}(\partial_-A)$,
\\
(2) $\gamma(1)\in \pi^{-1}(X\cap C)$, and
\\
(3) $\gamma([0,1))\subset \pi^{-1}(U_-\cap C)$.

Let $Y$ be an element of $\mathcal Y$ which contains $\gamma(1)$.
Then there are two unbounded
connected components of the  complement
$\pi^{-1}(C)\setminus(Y\cup\gamma)$,
one $L(Y\cup\gamma)$ on the left, and the other $R(Y\cup\gamma)$ on the right.

Notice that for any $n>0$, $\tilde g^{-n}\gamma$ is a path in $C$,
and that $p\tilde g^{-n}(\gamma(1))\to-\infty$ and
$p\tilde g^{-n}(\gamma(0))\to\infty$ as $n\to\infty$.
That is, for any large $n$, $\tilde g^{-n}(\gamma(1))\in L(Y\cup\gamma)$
and $\tilde g^{-n}(\gamma(0))\in R(Y\cup\gamma)$, showing that
$\tilde g^{-n}(\gamma)$
intersects $\gamma$. On the other hand in $\check U_-$,
$\gamma$ defines a curve from a point in $\pi^{-1}(\partial_-A)$
to a prime end in $\pi^{-1}(\partial_\infty U_-)$. But by
condition (4) of Lemma \ref{l2}, $\gamma$ cannot intersect
$\tilde g^{-n}(\gamma)$ for any large $n$.
A contradiction.

Finally the case where $\partial_-R$ is oriented from the left to the
right can be dealt with similarly by reversing the time.
This completes the proof of Theorem \ref{tx}.


\begin{thebibliography}{99} 

\bibitem[B]{B} G. D. Birkhoff,
{\em Sur quelques courbes ferm\'ees remarquables,} Bull.\ Soc.\
Math.\ France {\bf 60}(1932) 1-26; also in {\em Collected Mathematical
 Papers
of G. D. Birkhof,} vol.\ II, pp.\ 444-461

\bibitem[BG]{BG} M. Barge and R. M. Gillete, {\em Rotation and
periodicity in plane separating continua,} Ergod.\ Th.\ Dyn.\ Sys.\
{\bf 11}(1991) 619-631.

\bibitem[C]{C} M. Charpentier, {\em Sur quelques propri\'et\'es des courbes
de M. Birkhoff,} Bull.\ Soc.\ Math.\ France  {\bf 62}(1934) 193-224.



\bibitem[CL]{CL} M. L. Cartwright and J. E. Littlewood,
{\em Some fixed point theorems,} Ann.\ Math.\ {\bf 54}(1951)
1-37.

\bibitem[E]{E} D. B. A. Epstein, {\em Prime ends,} Proc.\ London
Math.\ Soc.\ {\bf 42}(1981) 385--414.

\bibitem[F]{F} J. Franks, {\em Recurrence and fixed points of surface 
homeomorphisms,} Ergod.\ Th.\ Dyn.\ Sys.\ {\bf 8}(1988) 99-107.

\bibitem[FL]{FL} J. Franks and P. Le Calvez, {\em Regions of
instability for non-twist maps,} Ergod.\ Th.\ Dyn.\ Sys.\ 
{\bf 23}(2003), 111--141.

\bibitem[L1]{L1} P. Le Calvez, {\em Existence d'orbits quasi-periodiques
dans les attracteurs de
Birkhoff,} Commun.\ Math.\ Phys.\ {\bf 106}(1986) 383-39.

\bibitem[L2]{L2} P. Le Calvez, {\em Propri\'et\'es des attracteurs
de Birkhoff,} Ergod.\ Th.\ Dyn.\ Sys.\ {\bf 8}(1987) 241-310

\bibitem[L3]{L3} P. Le Calvez, {\em Une version feuillet\'ee 
\'equivariante du th\'eor\`eme de translation de Brouwer,}
Publ.\ Math.\ I. H. E. S. {\bf 102}(2005) 1--98.



\bibitem[M]{M} J. Mather, {\em Topological proofs of some purely
topological consequences of Carath\'eodory's theory of prime ends,}
In: {Th. M. Rassias, G. M. Rassias, eds., Selected Studies,}
North-Holland, (1982) 225--255.



\bibitem[MN]{MN} S. Matsumoto and H. Nakayama, {\em Continua
as minimal sets of homeomorphisms of $S^2$,} Preprints in Arxiv.

     
\end{thebibliography}
\end{document}